\documentclass[aos,MSNbibl,nameyear,dvips]{arximspdf}

\doi{10.1214/13-AOS1111} 
\volume{41}
\issue{2}
\pubyear{2013}
\firstpage{1035}
\lastpage{1053}

\makeatletter
\def\cal{\mathcal}
\newtheorem{thmm}{Theorem}
\newtheorem{Lemma}{Lemma}
\newproclaim{remark}{Remark}
\makeatother

\begin{document}
\begin{frontmatter}

\title{Relative errors for bootstrap approximations of the serial
correlation coefficient}
\runtitle{Bootstrap approximations for serial correlation}

\begin{aug}
\author[A]{\fnms{Chris} \snm{Field}\ead[label=e1]{field@mathstat.dal.ca}\thanksref{e1}}
\and
\author[B]{\fnms{John} \snm{Robinson}\corref{}\ead[label=e2]{john.robinson@sydney.edu.au}\thanksref{e2}}
\thankstext{e1}{Supported by NSERC Discovery grant.}
\thankstext{e2}{Supported by ARC DP0773345.}
\runauthor{C. Field and J. Robinson}
\affiliation{Dalhousie University and University of Sydney}
\address[A]{Department of Mathematics and Statistics\\
Dalhousie University\\
Halifax, Nova Scotia\\
Canada B3H 3J5\\
\printead{e1}}
\address[B]{School of Mathematics and Statistics\\
University of Sydney\\
NSW 2006\\
Australia\\
\printead{e2}}
\end{aug}

\received{\smonth{3} \syear{2012}}
\revised{\smonth{11} \syear{2012}}

%
\begin{abstract}
We consider the first serial correlation coefficient under an $\operatorname{AR}(1)$
model where
errors are not assumed to be Gaussian. In this case it is necessary to consider
bootstrap approximations for tests based on the statistic since the distribution
of errors is unknown. We obtain saddle-point approximations for tail
probabilities
of the statistic and its bootstrap version and use these to show that
the bootstrap
tail probabilities approximate the true values with given relative
errors, thus
extending the classical results of Daniels [\textit{Biometrika} \textbf
{43} (1956) 169--185]
for the Gaussian case. The methods require conditioning on the set of
odd numbered
observations and suggest a conditional bootstrap which we show has
similar relative error properties.
\end{abstract}

%


\end{frontmatter}

\section{Introduction}\label{sec1}
A central limit theorem for the first-order serial correlation for an
autoregression with general errors was obtained by \citet{A59}, %
and Edgeworth expansions were obtained by \citet{BA}
who used this to prove the validity of the bootstrap approximation.
There have been several papers which consider saddle-point
approximations for
autoregressive processes [\citet{D56}, \citet{P78}, \citet{L94}] under the
assumption of normal errors and more generally for a ratio of
quadratic forms of normal variables [\citet{L89}]. Our results, in
contrast, give relative errors, valid for nonnormal errors and are
used to show that the bootstrap has better than first-order relative
accuracy in a moderately large region.

Let $\varepsilon_0,\varepsilon_1,\ldots,\varepsilon_n$ be independent and
identically distributed random variables with distribution function $F$
and density $f$, assume that $E\varepsilon_0=0$, define $X_i=\rho
X_{i-1}+\varepsilon_i, i=2,\ldots,n$ and take $X_1$ to be distributed as
$\varepsilon_0/\sqrt{1-\rho^2}$, which, although not of the correct form
of the stationary distribution when we do not assume normal errors, has
a variance in common with that case. We consider approximating the
distribution of the first serial correlation coefficient,
%
\begin{equation}
R=\frac{\sum_{i=2}^n X_iX_{i-1}}{X_1^2/2+\sum_{i=2}^{n-1}
X_i^2+X_n^2/2},\label{defR0}
\end{equation}
following Section 6 of \citet{D56} who obtained a saddle-point
approximation for this when $f$ was the density of a normal variable.
Note that without loss of generality we can assume $E\varepsilon_0^2=1$.
We wish to consider testing the hypothesis $\rho \leq\rho_0$ using $R$.

When $F$ is unknown we will consider a bootstrap approximation to the
test, generating a bootstrap sample, $X_1^*,\ldots,X_n^*$, under the
hypothesis using methods described later.
Then we can obtain $R^*$ by replacing $X_1,\ldots,X_n$ by $X_1^*,\ldots,X_n^*$ in the definition of $R$. We use a test based on $R^*$, so we
need to know the accuracy of the approximations $P^*(R^*>u)$ to $
P(R>u)$, where $P^*$ refers to probabilities under the bootstrap
sampling given the original sample.

We are unable to obtain a saddle-point approximation to this tail area
directly. Instead we will consider conditioning over a subset of the
random variables and obtain an approximation to the conditional tail
area. In order to get the unconditional tail area, we take the expected
value over the conditioning variables.
We will show that we can approximate the conditional distribution with a
saddle-point approximation where the conditioning is on~$\mathbf{C}$,
the odd
numbered observations. The approximation is
%
\begin{equation}
P(R \ge u|\mathbf{C})=\bar\Phi\bigl(\sqrt{m}W^+(u)\bigr) \bigl(1+O_P(1/m)
\bigr),\label{sadapp}
\end{equation}
where $m$ is the number of even numbered observations,$\bar\Phi
(z)=P(Z\geq z)$ for $Z$ a standard normal variable, and $W^+(u)$ is defined
later. We obtain a similar approximation for $P^*(R^*\geq u|\mathbf{C^*})$.

We want the relative error of the unconditional bootstrap tail area
under $\rho_0$ as an approximation of the true tail area. We use the
saddle-point
approximation as a device to enable this comparison. Since we cannot
get a saddle-point for the unconditional probability, we need
to work from the conditional approximations. Now $P(R\ge u)=EP(R\ge
u|\mathbf{C})$ and $P^*(R^*\ge u)=E^*P^*(R^*\ge u|\mathbf{C^*})$, where
$E^*$ is expectation under the bootstrap resampling given the original
sample. Then the relative error is
%
\begin{equation}
\frac{P(R \ge u)-P^*(R^* \ge u)}{P(R \ge u)}.\label{relerr}
\end{equation}

The above conditioning suggests a different conditional bootstrap, in
which we condition on the odd numbered observations $\mathbf{C}$ and
obtain conditional bootstrap samples for the even observations. This
permits a direct comparison of the conditional\vadjust{\goodbreak} distributions of the
ratios $R$ and a bootstrap counterpart given the same odd numbered
observations, $\mathbf{C}$. We describe this conditional bootstrap and
compare tests based on it to tests based on the unconditional
bootstrap. We introduce this conditional bootstrap and obtain a
saddle-point approximation for it.

The next section provides the details of the conditioning and is
followed by a section giving results for the Gaussian case for both
conditional and unconditional cases, then by sections giving the
derivation of the main result. A final section provides some numerical
results illustrating the accuracy of the approximations and comparing
the power of the conditional and unconditional bootstraps.

\section{Conditioning}\label{seccond}

Assume that $n=2m+1$. Let
\[
S=\sum_{i=2}^n X_iX_{i-1}-u
\Biggl(X_1^2/2+ \sum_{i=2}^{n-1}
X_i^2+X_n^2/2\Biggr),
\]
then $P(R>u)=P(S>0)$.
Let
$A_i=X_{2i-1}+X_{2i+1}$, $B_i=(X_{2i-1}^2+X_{2i+1}^2)/2$ for $i=1,\ldots,m$,
and $\mathbf{C}=(X_1,X_3,\ldots,X_n)$, and write
%
\begin{eqnarray}
\label{SC} S&=&\sum_{i=1}^{m}
\bigl(A_i X_{2i}-u\bigl(X_{2i}^2+B_i
\bigr)\bigr)
\nonumber
\\[-8pt]
\\[-8pt]
\nonumber
&=&-u\sum_{i=1}^{m} (X_{2i}-A_i/2u)^2+
m\frac{\bar{A}^2-4u^2\bar
B}{4u},
\end{eqnarray}
where $m\bar{A}^2=\sum_{i=1}^{m}A_i^2$ and $m\bar B=\sum_{i=1}^{m}B_i$.
So for $u>0$,
$P(S>0|\mathbf{C})=0$ if $\bar{A}^2-4u^2\bar B<0$.

It is clear that when $\rho_0=0$, conditional on $\mathbf{C}$, the
terms in the
sums in $S$ are independent random variables.
If $\rho_0\ne0$ the first step is to show that the $X_{2i}$'s are independent
conditional on $\mathbf{C}$.
This follows since we can factor the joint density of $\mathbf
{D}=(X_2,X_4, \ldots, X_{n-1})$
conditional on $\mathbf{C}=(X_1,X_3,\ldots,X_n)$.

\section{The Gaussian case}\label{sec3}
We will first give a brief account of the saddle-point approximations
for the Gaussian case where both an unconditional and conditional
approach are possible with explicit forms for the approximations.

Consider the unconditional normal case. If $\varepsilon_1,\ldots,\varepsilon
_n$ are independent standard normal, $X_1=\varepsilon_1/\sqrt{1-\rho^2}$
and $X_i=\rho X_{i-1}+\varepsilon_i$ for $i=2,\ldots,n$, and
\[
S=\sum_{i=2}^n X_iX_{i-1}-u
\Biggl(X_1^2/2+\sum_{i=2}^{n-1}
X_i^2+X_n^2/2
\Biggr)=x^T(A-uB)x,
\]
with $A$ and $B$ symmetric. We find the saddle-point approximation to
\mbox{$P(S\geq0)$} following the method of \citet{L94}. The cumulative
generating function of $S$ is
\begin{eqnarray*}
\kappa(t)&=&\log\bigl((2\pi)^{n/2}|\Sigma|^{1/2}
\bigr)^{-1}\int e^{tx^T(A-uB)x-x^T\Sigma^{-1}x/2}\,dx
\\
&=&\log\bigl|I-2tU(A-uB)U^T\bigr|^{-1/2}
\\
&=&-\frac{1}{2}\sum_{i=1}^n
\log(1-2t\lambda_i),
\end{eqnarray*}
where $\sigma_{ij}=\rho^{|i-j|}$, $\Sigma=U^TU$, $U$ is upper
triangular and $\lambda_1\leq\cdots\leq\lambda_n$ are the eigenvalues
of $U(A-uB)U^T$. So the Barndorff--Nielsen approximation [see
Section 1.2 of \citet{FRS13}] is
\[
P(S\geq0)=\bar\Phi\bigl(\sqrt{m}w^\dagger\bigr) \bigl(1+O(1/n)\bigr),
\]
where $w^\dagger=w-\log\psi(w)/nw$ for $w=(-2\kappa(\hat t))^{1/2}$,
where $\hat t$ is the solution to $\kappa'(\hat{t})=0$ and $\psi(w)=w/\hat
t(\kappa''(\hat t))^{1/2}$. Note that $\kappa(t)$, $\hat t$, $w$ and so
$w^\dagger$ all are functions of $u$, but this dependence is suppressed
to simplify notation.

To consider the power of the test $H_0\dvtx \rho=\rho_0$ versus the
alternative $H_1\dvtx \rho=\rho_1>\rho_0$, we can
find the critical values from the saddle-point approximation under
$H_0$ for a fixed level and then the power directly under $H_1$.

Now consider the conditional test. If the observations are as above and
$A_1,\ldots,A_m$ and $B_1,\ldots,B_m$ are defined as in Section \ref
{seccond}, then we need to find $P(S\geq0|\mathbf{C})$. Recall that
\[
S=\sum_{i=1}^m\bigl(X_{2i}A_i-u
\bigl(X_{2i}^2+B_i\bigr)\bigr),
\]
and in this case, given $A_i$ and $B_i$, $X_{2i}$ are conditionally
independent with conditional distribution normal with mean $\rho
A_i/(1+\rho^2)$ and variance\break $1/(1+\rho^2)$.
The test of $H_0$ will be performed by considering the conditional
distribution of $S$ given $\mathbf{C}$ obtained when $X_{2i}$ are
assumed to be conditionally independent normal variables with mean $\rho
_0A_i/(1+\rho_0^2)$ and variance $1/(1+\rho_0^2)$. So the critical
value at a fixed level can be calculated from this distribution. Then
the power can be calculated using the conditional distribution of $S$
given $\mathbf{C}$ using $X_{2i}$ conditionally independent normal
variables with mean $\rho_1 A_i/(1+\rho_1^2)$ and variance $1/(1+\rho
_1^2)$. These conditional distributions can be approximated by a
saddle-point method as in the unconditional case, by using the
conditional cumulative generating function of $S$, given by
%
\begin{eqnarray}
\label{knorcond} \kappa(t)&=&\frac{1}{m}\sum_{i=1}^m
\log\sqrt{\frac{1+\rho^2}{2\pi}}\int e^{-tu(z-A_i/2u)^2-(1+\rho^2)(z-\rho A_i/(1+\rho^2))^2/2}\,dz\nonumber\\
&&{} +t\frac{\bar{A}^2-4u^2\bar
{B}}{4u}
\\
&=&-\frac{1}{2}\log\biggl(1+\frac{2tu}{1+\rho^2}\biggr)-tu\bar{B}+
\frac{\bar
{A}^2(\rho+t)^2}{2(1+\rho^2+2tu)}-\frac{\bar{A}^2\rho^2}{2(1+\rho
^2)}.\nonumber
\end{eqnarray}
From (\ref{knorcond}), $\kappa(0)=0$, and differentiating (\ref
{knorcond}) shows that for $u>0$, $\kappa'(0)<0$ and that $\kappa
'(t)<0$ for all $t>0$ if $\bar{A}^2-4u^2\bar{B}<0$ and that $\kappa
'(t)\to(\bar{A}^2-4u^2\bar{B})/4u$ as $t\to\infty$. So $\kappa'(t)=0$
has a solution, if and only if $\bar{A}^2-4u^2\bar B >0$. Then the
Barndorff--Nielsen approximation for the conditional distribution can
be obtained as before.

\section{The general case}\label{sec4}

We can get a general bootstrap sample by considering the residuals
$\varepsilon_i=X_i-\rho_0X_{i-1}$, $i=2,\ldots,n$ and drawing bootstrap
replicates by sampling $\varepsilon_1^*,\ldots,\varepsilon_n^*$ from
$F_n(x)=\sum_{i=2}^nI((\varepsilon_i-\bar\varepsilon)/\sigma_n\leq x)/(n-1)$,
where $ \bar\varepsilon=\sum_{i=2}^n\varepsilon_i/(n-1)$ and $\sigma_n^2=\sum_{i=2}^n(\varepsilon_i-\bar\varepsilon)^2/(n-1)$, then generating bootstrap
versions of the sample as $X_1^*=\varepsilon_1^*/\sqrt{1-\rho_0^2}$,
$X_i^*=\rho_0X_{i-1}^*+\varepsilon_i^*$ for $i=2,\ldots,n$. From this
bootstrap sample we can calculate $R^*$ unconditionally.\looseness=1

We consider saddle-point approximations to the conditional distribution
of $S$ given $\mathbf{C}$ then get the approximation to the
unconditional distribution by considering the expectation of these. For
the bootstrap no density exists, so we consider a smoothed bootstrap by
adding independent normal variables with zero mean and small standard
deviation $\tau$ to each bootstrap value $\varepsilon_1^*,\ldots,\varepsilon
_n^*$ obtaining $\varepsilon_1^\dagger,\ldots,\varepsilon_n^\dagger$. Then we
can proceed in the same way to approximate the bootstrap distribution
as the expectation of the approximation to the conditional
distribution. Finally we show that for a suitable choice of $\tau$ the
smoothed bootstrap approximates the unconditional bootstrap with
appropriate relative error.

We also consider a conditional bootstrap where we condition on $\mathbf
{C}$, the same conditioning variables used for the true distribution.
Here we are able to obtain relative errors for the approximation to the
conditional distribution of $S$ given $\mathbf{C}$.

\subsection{Approximations under conditioning}\label{appcond}
From the factorization
of the joint density of $\mathbf{D}=(X_2,X_4, \ldots,X_{n-1})$
conditional on $\mathbf{C}=(X_1,X_3,\ldots,X_n)$,
we get the conditional
density of $X_{2i}$ given $X_{2i-1}$ and $X_{2i+1}$ is
\begin{eqnarray*}
&&g(z|X_{2i-1},X_{2i+1})\\
&&\qquad= f(z|X_{2i-1})f(X_{2i+1}|z)/f(X_{2i+1}|X_{2i-1})
\\
&&\qquad=\frac{f_{\varepsilon}(z-\rho_0 X_{2i-1})
f_{\varepsilon}(X_{2i+1}-\rho_0 z)}{\int f_{\varepsilon}(z-\rho_0 X_{2i-1})
f_{\varepsilon}(X_{2i+1}-\rho_0 z)\,dz},
\end{eqnarray*}
where $f_{\varepsilon}$ is the density of the errors $\varepsilon_2,\ldots,\varepsilon_n$. Define $S$ as in (\ref{SC}). Then we can get
approximations to the distribution of $S$ given $\mathbf{C}$ using this density.

The conditional cumulant generating function for $S$ given $\mathbf{C}$ is
%
\begin{eqnarray}
\label{Kr} mK(t,u)&=&\sum_{i=1}^{m}\log \int
e^{\{t(A_i z-u(z^2+
B_i))\}}g(z|X_{2i-1},X_{2i+1})\,dz
\nonumber\\
&=&\sum_{i=1}^{m}\log \int
e^{-tu(z-A_i/2u)^2}g(z|X_{2i-1},X_{2i+1})\,dz \\
&&{}+m\frac{t(\bar{A}^2-4u^2\bar B)}{4u}.\nonumber
\end{eqnarray}
Note that this will exist whenever $tu>0$.
We use the notation $K_{ij}(t,u)=\partial^{i+j}K(t,u)/\partial
t^i\partial u^j$.
Then differentiating (\ref{Kr}) with respect to $t$ gives
%
\begin{equation}
K_{10}(t,u)=-\frac{1}{m}\sum_{i=1}^{m}
K_i(t,u)+\frac{(\bar
{A}^2-4u^2\bar B)}{4u}\label{dK}
\end{equation}
and
%
\begin{eqnarray}
\label{K20tu} K_{20}(t,u)&=&\frac{1}{m} \sum
_{i=1}^{m}\frac{\int
u^2(z-A_i/2u)^4e^{-tu(z-A_i/2u)^2}g(z|X_{2i-1},X_{2i+1})\,dz}{\int
e^{-tu(z-A_i/2u)^2}g(z|X_{2i-1},X_{2i+1})\,dz}
\nonumber
\\[-8pt]
\\[-8pt]
\nonumber
&&{}-\frac{1}{m} \sum_{i=1}^{m}K_i(t,u)^2,
\end{eqnarray}
where
%
\begin{equation}
K_i(t,u)=\frac{\int
u(z-A_i/2u)^2e^{-tu(z-A_i/2u)^2}g(z|X_{2i-1},X_{2i+1})\,dz}{\int
e^{-tu(z-A_i/2u)^2}g(z|X_{2i-1},X_{2i+1})\,dz}.\label{Ki}
\end{equation}
Note from (\ref{Kr}) that $K(0,u)=0$ and from (\ref{dK}) that if $\bar
{A}^2-4u^2\bar B<0$, then $K_{10}(t,u)$ is always negative, so there is
no solution to the saddle-point equation \mbox{$K_{10}(t,u)=0$}. For $\bar
{A}^2-4u^2\bar B>0$ we first find a value of $u$ such that
$K_{10}(0,u)=0$. Now
\[
K_{10}(0,u)=\frac{1}{m} \sum_{i=1}^{m}
\int \bigl(zA_i-uz^2\bigr)g(z|X_{2i-1},X_{2i+1})\,dz-u
\bar B.
\]
Let $u_{0}$ be such that $K_{10}(0,u_{0})=0$, then
%
\begin{equation}
u_0=\frac{ \sum_{i=1}^{m}\int zg(x|X_{2i-1},X_{2i+1})\,dz A_i}{\sum_{i=1}^{m}\int z^2g(z|X_{2i-1},X_{2i+1})\,dz+m\bar B}.\label{u0}
\end{equation}
So for $u> u_0$,
\[
K_{10}(0,u)=(u_0-u) \Biggl(\frac{1}{m}\sum
_{i=1}^{m}\int z^2g(z|X_{2i-1},X_{2i+1})\,dz+
\bar B\Biggr)<0
\]
and $K_{20}(t,u)>0$. So for $u>u_0$, $K_{10}(t,u)$ is increasing in
$t$, is negative for $t=0$ and as $t\to\infty$,
\[
K_{10}(t,u)\to\frac{\bar{A}^2-4u^2\bar B}{4u},
\]
since the first term in (\ref{dK}) tends to 0 as $t\to\infty$. Thus the
saddle-point equation $K_{10}(t,u)=0$, has a finite solution, $t(u)$
for $u>u_0$, if and only if \mbox{$\bar{A}^2-4u^2\bar B>0$}. Further,
$K(t(u),u)$ exists and is finite if $\bar{A}^2-4u^2\bar B>0$. If $\bar
{A}^2-4u^2\bar B<0$, $K(t,u)\to-\infty$ as $t\to\infty$.

If $\bar{A}^2-4u^2\bar B>0$, the Barndorff--Nielsen form of the
saddle-point approximation is
%
\begin{equation}
P (S\geq0 |\mathbf{C}=\mathbf{c} ) =\bar\Phi\bigl(\sqrt{m} W^+\bigr)
\bigl(1+O_P\bigl(m^{-1}\bigr)\bigr),\label{conpr}
\end{equation}
where
%
\begin{equation}
W^+=W-\log\bigl(\Psi(W)\bigr)/(m W),\label{W+}
\end{equation}
with
%
\begin{equation}
W=\sqrt{-2K\bigl(t(u),u\bigr)}\quad\mbox{and} \quad\Psi(W)=W/\bigl(t(u)
\sqrt{K_{20}\bigl(t(u),u\bigr)}\bigr).\label{psi0}
\end{equation}
The proof of this result is given in Section 1 of the supplementary
material of \citet{FRS13}.

The bootstrap distribution of $\varepsilon_1^*,\ldots,\varepsilon_n^*$ does
not have a density, but we can approximate the distribution by a
smoothed version which is continuous. Let
%
\begin{equation}
f_n(z)=\frac{1}{n-1}\sum_{k=2}^n
\frac{e^{-(z-\eta_k)^2/2\tau^2}}{\sqrt {2\pi\tau^2}},\label{fn}
\end{equation}
where $\eta_k=(\varepsilon_k-\bar\varepsilon)/\sigma_n$.
If we draw a sample $\varepsilon^\dagger_1,\ldots,\varepsilon^\dagger_n$ from
this distribution and obtain
$X_1^\dagger=\varepsilon^\dagger_1/(1-\rho_0^2)$ and $X_i^\dagger=\rho_0
X_{i-1}^\dagger+\varepsilon^\dagger_i$,
then choosing $\tau$ small enough, we can approximate the bootstrap
distribution of $R^*$ by the bootstrap version of~$R^\dagger$. With
this new smoothed bootstrap we can proceed to get the saddle-point
approximation to its distribution by using the expectation of the
conditional bootstrap as we do for the saddle-point approximation of
the distribution of $R$.

The conditional density of $X_{2i}^\dagger$ given $X_{2i-1}^\dagger$
and $X_{2i+1}^\dagger$ is
%
\begin{equation}
g^\dagger\bigl(z|X_{2i-1}^\dagger,X_{2i+1}^\dagger
\bigr)=\frac{f_n(z-\rho_0
X_{2i-1}^\dagger)f_n(X_{2i+1}^\dagger-\rho_0 z)}{\int f_n(z-\rho_0
X_{2i-1}^\dagger)f_n(X_{2i+1}^\dagger-\rho_0 z)\,dz},\label{gdag}
\end{equation}
where
%
\begin{equation}
g^\dagger\bigl(z|X_{2i-1}^\dagger,X_{2i+1}^\dagger
\bigr)=\frac{1}{(n-1)^2}\sum_k\sum
_l g_{ikl}^\dagger(z)\label{gidag}
\end{equation}
for
%
\begin{equation}
g_{ikl}^\dagger(z)=\frac{(n-1)^2e^{-(z-\rho_0 X_{2i-1}^\dagger-\eta
_k)^2/2\tau^2-(X_{2i+1}^\dagger-\rho_0 z-\eta_l)^2/2\tau^2}}{\sum_k\sum_l\int e^{-(z-\rho_0 X_{2i-1}^\dagger-\eta_k)^2/2\tau
^2-(X_{2i+1}^\dagger-\rho_0 z-\eta_l)^2/2\tau^2}\,dz}.\label{gikldag}
\end{equation}
Now
\begin{eqnarray*}
&&\bigl[\bigl(z-\rho_0 X_{2i-1}^\dagger-
\eta_k\bigr)^2+\bigl(X_{2i+1}^\dagger-
\rho_0 z-\eta _l\bigr)^2\bigr]
\\
&&\qquad=\bigl(1+\rho_0^2\bigr) \biggl(z'-
\frac{\eta_k-\rho_0\eta_l }{1+\rho_0^2}\biggr)^2\\
&&\qquad\quad{}+\frac
{(X_{2i+1}^\dagger-\rho_0^2 X_{2i-1}^\dagger-\rho_0\eta_k-\eta
_l)^2}{(1+\rho_0^2)},
\end{eqnarray*}
where $z'=z-\rho_0(X_{2i-1}^\dagger+X_{2i+1}^\dagger)/(1+\rho_0^2)$.
So, integrating with respect to $z$ in the denominator of
$g_{ikl}^\dagger(z)$ we have
\begin{eqnarray*}
g_{ikl}^\dagger(z)&=&\frac{e^{-(1+\rho_0^2)(z'-(\eta_k-\rho_0\eta
_l)/(1+\rho_0^2))^2/2\tau^2-(X_{2i+1}^\dagger-\rho_0^2 X_{2i-1}^\dagger
-\rho_0\eta_k-\eta_l)^2/2\tau^2(1+\rho_0^2)}}{\sqrt{2\pi\tau^2}\sum_k\sum_l e^{-(X_{2i+1}^\dagger-\rho_0^2 X_{2i-1}^\dagger-\rho_0\eta
_k-\eta_l)^2/2\tau^2(1+\rho_0^2)}}.
\end{eqnarray*}
Define $S^\dagger$ as in (\ref{SC}) using $X^\dagger$ in place of $X$,
with analogous definitions for $A_i^\dagger$, $B_i^\dagger$, $R^\dagger
$ and $C^\dagger$. Then the conditional cumulant generating function of
$S^\dagger$ given $\mathbf{C}^\dagger$ is
%
\begin{eqnarray}
\label{Kdag} mK^\dagger(t,u) &=&\sum_{i=1}^{m}\log
\int e^{-tu(z-A_i^\dagger/2u)^2}g^\dagger \bigl(z|X_{2i-1}^\dagger,X_{2i+1}^\dagger
\bigr)\,dz
\nonumber
\\[-8pt]
\\[-8pt]
\nonumber
&&{} + m\frac{t(\bar{A}^{\dagger2}-4u^2\bar{B^\dagger
})}{4u},
\end{eqnarray}
which is of the same form as the formula for $K(t,u)$ with $g^\dagger
(z|X_{2i-1}^\dagger,X_{2i+1}^\dagger)$ replacing
$g(z|X_{2i-1},X_{2i+1})$. So we can obtain analogous results to those
of (\ref{dK})--(\ref{u0}) and to the argument following these, to show
that, when $\bar{A}^{\dagger2}-4u^2\bar B^\dagger>0$, if $t^\dagger
(u)$ is the solution of
$K_{10}^\dagger(t,u)=0$, then the saddle-point approximation is
\[
P^\dagger\bigl(S^\dagger\geq0 |\mathbf{C}^\dagger\bigr) =
\bar\Phi\bigl(\sqrt{m} W^{\dagger+}\bigr) \bigl(1+O_P
\bigl(m^{-1}\bigr)\bigr),
\]
where
\[
W^{\dagger+}=W^*-\log\bigl(\Psi^\dagger(W)\bigr)/\bigl(m
W^\dagger\bigr),
\]
with
\[
W^\dagger=\sqrt{-2K^\dagger\bigl(t^\dagger(u),u\bigr)}
\quad\mbox{and}\quad \Psi^\dagger(W)=W^\dagger/\bigl(t^*(u)\sqrt
{K_{20}^\dagger\bigl(t^\dagger(u),u\bigr)}\bigr).
\]
We can summarize these results in the following theorem:

\begin{thmm}\label{rne0condthm}
For $u\geq0$, $P(S>0|\mathbf{C})=0$ if $\bar{A}^{2}-4u^2\bar B<0$ and
$P(S^\dagger>0|\mathbf{C}^\dagger)=0$ if $\bar{A}^{\dagger2}-4u^2\bar
B^\dagger<0$. If
$\bar{A}^{2}-4u^2\bar B>0$ and $\bar{A}^{\dagger2}-4u^2\bar B^\dagger
>0$, for $u>u_0$ from (\ref{u0}) and $u>u_{0}^\dagger$ defined analogously,
$t(u)$ and $t^\dagger(u)$, solutions of $K_{10}(t,u)=0$ and
$K_{10}^\dagger(t,u)=0$,
exist and are both finite and positive, and if $EX_1^8$ is bounded,
\[
P(R >u|\mathbf{C})= \bar\Phi\bigl(\sqrt{m}W^+\bigr)\bigl[1+O_P(1/m)
\bigr]
\]
and
\[
P^\dagger\bigl(R^\dagger>u|\mathbf{C}^\dagger\bigr)= \bar
\Phi\bigl(\sqrt{m}W^{\dagger
+}\bigr)\bigl[1+O_P(1/m)\bigr],
\]
where $W(u)$, $W^+$ and $\Psi(W^\dagger)$ are defined as in (\ref{W+})
and (\ref{psi0}) and
\[
W^{\dagger+}=W^\dagger-\log\bigl(\Psi^\dagger\bigr)/\bigl(m
W^\dagger\bigr),
\]
with
\[
W^\dagger=\sqrt{-2K^\dagger\bigl(t^\dagger(u),u\bigr)}
\quad\mbox{and}\quad \Psi^\dagger\bigl(W^\dagger\bigr)=W^\dagger/
\bigl(t^\dagger(u)\sqrt{K_{20}^\dagger
\bigl(t(u),u\bigr)}\bigr).
\]
\end{thmm}

\begin{remark*}
If $R'$ has the denominator in $R$ replaced by $\sum_{i=1}^nX_i^2$, then $P(R'>u|\mathbf{C})=P(S>u(X_1^2+X_n^2)/2|\mathbf
{C})$. So we can proceed with the saddle-point approximation obtaining
results with the relative error unchanged, since throughout the errors
will be affected by a term of $O_P(u/m)$. A similar argument gives
results for $n$ even.
\end{remark*}

\subsection{The relative error of the bootstrap}\label{relerrbsnonzero}
Assume throughout this section that the conditions of Theorem~\ref{rne0condthm} hold.
Let ${{\cal{A}}}=\{\mathbf{C}\dvtx \bar{A}^{2}-4u^2\bar B>0\}$.
Now $E(\bar{A}^2-4u^2\bar B)=(2(1-2u^2)+2\rho_0^2)/(1-\rho_0^2)$ and
$\operatorname{var}(\bar{A}^2-4u^2\bar B)=O(1/m)$, so for $1-2u^2+\rho
_0^2>\delta>0$, it follows from the Chebychev inequality that
$P({\cal{A}}^c)=P(\bar{A}^2-4u^2\bar B<0)=O(1/m)$. So, since
$P(S>0|\mathbf{C})I({\cal{A}}^c)=0$,
%
\begin{eqnarray}
\label{PRgu} P(S>0) &=&E\bigl[P(S>0|\mathbf{C})I({\cal{A}})\bigr]+E\bigl[P(S>0|
\mathbf{C})I\bigl({\cal {A}}^c\bigr)\bigr]
\nonumber
\\[-8pt]
\\[-8pt]
\nonumber
&=&E\bigl[\bar\Phi\bigl(\sqrt{m}W^+\bigr)I({\cal{A}}) \bigl(1+O_P(1/m)
\bigr)\bigr].
\end{eqnarray}
Restrict attention to $\cal{A}$, so with $u_0$ given in (\ref{u0}),
$K_{10}(0,u_0)=0$ and thus $t(u_0)=0$ and
\[
t(u)=t'(u_o) (u-u_0)+
\tfrac{1}{2}t''(u_0)
(u-u_o)^2+O_P\bigl((u-u_0)^3
\bigr).
\]
Further, since $K_{10}(t(u),u)=0$,
$t'(u_0)=-K_{11}/K_{20}$,
where we write $K_{ij}=K_{ij}(0,u_0)$. Then expanding $K(t(u),u)$ about
$u_0$ we obtain,
\[
K\bigl(t(u),u\bigr)=-D_1(u-u_0)^2-D_2(u-u_0)^3+O_P
\bigl((u-u_0)^4\bigr),
\]
where $D_1=K_{11}^2/2K_{20}$ and
%
\begin{equation}
D_2=\tfrac{1}{2}\bigl[t''_0K_{11}+t'_0K_{12}+t^{\prime2}_0K_{21}+
\tfrac
{1}{3}t^{\prime3}_0K_{30}\bigr].
\label{B20}
\end{equation}
So
%
\begin{equation}
W=(u-u_0)\sqrt{2D_1}\bigl(1+(u-u_0)D_2/2D_1
\bigr)+O_P\bigl((u-u_0)^3\bigr).
\label{Wexp}
\end{equation}

Note that $u_0$ is given in (\ref{u0}), so
%
\begin{equation}
u_0=\frac
{E[E(X_{2}|X_1,X_3)(X_1+X_3)]}{E[E(X_{2}^2|X_1,X_3)+X_1^2]}+J_u/\sqrt {m}+O_P(1/m),
\label{u0app}
\end{equation}
where, here and in the sequel, values of $J$ denote zero mean random
variables with finite variances.
Further, since $X_2=\rho_0X_1+\varepsilon_2$, $X_3=\rho_0^2X_1+\rho
_0\varepsilon_2+\varepsilon_3$ and $X_1$ is independent of $\varepsilon_2$ and
$\varepsilon_3$, the numerator in (\ref{u0app}) is
\[
\rho_0EX_1(X_1+X_3)+E
\bigl[E(\varepsilon_2|\rho_0\varepsilon_2+
\varepsilon_3) \bigl(\rho _0^2X_1+
\rho_0\varepsilon_2+\varepsilon_3\bigr)
\bigr],
\]
and since $\varepsilon_2=((\varepsilon_2-\rho_0\varepsilon_3)+\rho_0(\rho
_0\varepsilon_2+\varepsilon_3))/(1+\rho_0^2)$, the numerator is
\[
\frac{\rho_0}{1-\rho_0^2}+\frac{\rho_0^3}{1-\rho_0^2}+\rho_0=\frac{2\rho
_0}{1-\rho_0^2}.
\]
The denominator of (\ref{u0app}) is
\[
E\bigl[E\bigl(X_{2}^2|X_1,X_3
\bigr)+X_1^2\bigr]=E\bigl(X_2^2
\bigr)+E\bigl(X_1^2\bigr)=\frac{2}{1-\rho_0^2}.
\]
So
%
\begin{equation}
u_0=\rho_0+J_u/\sqrt{m}+O(1/m).
\label{u0exp}
\end{equation}
From (\ref{dK}) and (\ref{Ki})
\[
K_{11}=-\frac{1}{m}\sum_{i=1}^m
\int z^2g(z|X_{2i-1},X_{2i+1})\,dz-\bar B
\]
so
\[
EK_{11}=-E\bigl(X_2^2+X_1^2
\bigr)=-\frac{2}{1-\rho_0^2}
\]
and
%
\begin{equation}
K_{11}=-\frac{2}{1-\rho_0^2}+J_{11}/\sqrt{m}+O_P(1/m).
\label{K11exp}
\end{equation}
From (\ref{K20tu}), and using (\ref{u0exp}), we can write
%
\begin{eqnarray}
\label{K20exp} K_{20}&=&\frac{1}{m}\sum
_{i=1}^m\biggl\{\int\bigl(\rho _0z^2-A_iz
\bigr)^2g(z|X_{2i-1},X_{2i+1})\,dz
\nonumber
\\
&&\hspace*{32pt}{} -\biggl[\int\bigl(\rho_0z^2-A_iz
\bigr)g(z|X_{2i-1},X_{2i+1})\,dz\biggr]^2\biggr
\}
\nonumber
\\[-8pt]
\\[-8pt]
\nonumber
&&{}+J_{20}/\sqrt {m}+O_P(1/m)
\\
&=&\frac{1}{m}\sum_{i=1}^m
\gamma(X_{2i-1},X_{2i+1})+J_{20}/\sqrt
{m}+O_P(1/m),\nonumber
\end{eqnarray}
so
%
\begin{eqnarray}
K_{20}&=&E_{20}+J_{20}'/
\sqrt{m}+O_P(1/m),\label{K20exp'}
\end{eqnarray}
where
\[
E_{20}=\frac{1}{m}\sum_{i=1}^m
E \gamma(X_{2i-1},X_{2i+1}).
\]
Now, recalling that $D_1=K_{11}^2/2K_{20}$, and using (\ref{K11exp})
and (\ref{K20exp'}), we have
%
\begin{equation}
D_1=\frac{2}{(1-\rho_0^2)^2E_{20}}+J_D/\sqrt{m}+O_P(1/m),
\label{D1exp}
\end{equation}
$t(u)=-(u-u_0)K_{11}/K_{20}+O_P((u-u_0)^2)$,
$\Psi(u)=W/t(u)\sqrt{K_{20}}=1+O_P(u-u_0)$, so
$\log\Psi(u)/mW=O_P(1/m)$,
and, from (\ref{W+}), (\ref{Wexp}), (\ref{u0exp}) and (\ref{D1exp}),
%
\begin{eqnarray}
\label{W+-EW+} &&W^+-EW^+
\nonumber
\\[-8pt]
\\[-8pt]
\nonumber
&&\qquad=(u-\rho_0) \biggl(\frac{J_W}{\sqrt{m}}+(u-
\rho_0)\frac
{H}{\sqrt{m}}\biggr)+O_P\biggl((u-
\rho_o)^3+\frac{1}{m}\biggr),
\end{eqnarray}
where $H=\sqrt{m}(D_2/2D_1-ED_2/2ED_1)$.

We can consider the smoothed bootstrap introduced in Section \ref
{appcond} in the same way. Let $W^\dagger$, $W^{\dagger+}$ be defined
as in the statement of Theorem~\ref{rne0condthm}, and let
$\cal{A}^\dagger=\{\mathbf{C}^\dagger\dvtx {\bar{A}^{\dagger2}-4u^2\bar
{B}^\dagger}>0\}$ and $E_+^\dagger(\cdot)=E^\dagger(\cdot|\cal{A}^\dagger)$.
Then restricting attention to $\cal{A}^\dagger$, $K_{10}^\dagger
(0,u_0^\dagger)=0$, so $t^{\dagger}(u_0^\dagger)=0$ and
\[
t^\dagger(u)=t^{\dagger\prime}\bigl(u_0^\dagger\bigr)
\bigl(u-u_0^\dagger\bigr)+\tfrac
{1}{2}t^{\dagger\prime\prime}
\bigl(u_0^\dagger\bigr) \bigl(u-u_0^\dagger
\bigr)^2+O_P\bigl(\bigl(u-u_0^\dagger
\bigr)^3\bigr),
\]
with
\[
t^{\dagger\prime}\bigl(u_0^\dagger\bigr)=-K_{11}^\dagger/K_{20}^\dagger,
\]
where $K_{ij}^\dagger=K_{ij}^\dagger(0,u_0).$
Now we proceed as above with $X_i^\dagger$, $g_i^\dagger
(\cdot|X_{2i-1}^\dagger,X_{2i+1}^\dagger))$, $E^\dagger(\cdot)$ and $E^\dagger
(\cdot|\cdot)$ replacing $X_i$, $g(z|X_{2i-1},X_{2i+1})$, $E(\cdot)$ and $E(\cdot|\cdot)$. So
\begin{eqnarray}
u_0^\dagger&=&\rho_0+J_u^\dagger/
\sqrt{m}+O_P\biggl(\frac{\rho_0}{\sqrt{m}}\biggr),
\nonumber\\
K_{11}^\dagger&=&-\frac{2}{1-\rho_0^2}+J_{11}^\dagger/
\sqrt{m}+O_P(1/\sqrt {m})\label{K11dagexp}
\end{eqnarray}
and
%
\begin{eqnarray}
\label{K20dagexp} K_{20}^\dagger&=&\frac{1}{m}\sum
_{i=1}^m\biggl\{\int\bigl(
\rho_0 z^2-A_i^\dagger z
\bigr)^2g^\dagger\bigl(z|X_{2i-1}^\dagger,X_{2i+1}^\dagger
\bigr)\,dz
\nonumber
\\
&& \hspace*{32pt}{}-\biggl[\int\bigl(\rho_0 z^2-A_i^\dagger
z\bigr)g^\dagger\bigl(z|X_{2i-1}^\dagger,X_{2i+1}^\dagger\bigr)\,dz\biggr]^2\biggr
\}
\nonumber
\\[-8pt]
\\[-8pt]
\nonumber
&&{}+J_{20}^\dagger/\sqrt{m}+O_P(1/m)
\\
&=& \frac{1}{m}\sum_{i=1}^m
\gamma^\dagger\bigl(X_{2i-1}^\dagger,X_{2i+1}^\dagger
\bigr)+J_{20}^\dagger/\sqrt{m}+O_P(1/m).\nonumber
\end{eqnarray}

In order to compare the first terms of (\ref{K20exp}) and (\ref
{K20dagexp}), we need first to replace $\gamma^\dagger(\cdot)$ in this
first term by $\gamma(\cdot)$ appearing in $E_{20}$. The following lemma,
the proof of which is given in Section 2 of the supplementary material of \citet
{FRS13}, accomplishes this.
%
\begin{Lemma}\label{le1} For $\tau=O(1/\sqrt{m})$,
\[
\int h(z)g^\dagger\bigl(z|X_{2i-1}^\dagger,X_{2i+1}^\dagger
\bigr)\,dz=\int h(z)g\bigl(z|X_{2i-1}^\dagger,X_{2i+1}^\dagger
\bigr)\,dz+\frac{J_h}{\sqrt {m}}+O_P\biggl(\frac{1}{m}\biggr).
\]
\end{Lemma}

Using Lemma~\ref{le1},
\[
\frac{1}{m}\sum_{i=1}^m
\gamma^\dag\bigl(X_{2i-1}^\dagger,X_{2i+1}^\dagger
\bigr)= \frac{1}{m}\sum_{i=1}^m
\gamma\bigl(X_{2i-1}^\dagger,X_{2i+1}^\dagger\bigr)
+\frac{J_h}{\sqrt{m}}+O_P\biggl(\frac{1}{m}\biggr),
\]
so
%
\begin{equation}
E_{20}^\dagger=\frac{1}{m}\sum
_{i=1}^mE\gamma^\dag\bigl(X_{2i-1}^\dagger,X_{2i+1}^\dagger\bigr)=E_{20}+J_{20}^\ddagger/
\sqrt{m}+O_P(1/m).\label{E20dagtoE20}
\end{equation}
Now, as before $D_1^\dagger=K_{11}^{\dagger2}/2K_{20}^\dagger$, so
using (\ref{K11dagexp}) and (\ref{E20dagtoE20}), we have
\[
D_1^\dagger=\frac{2}{(1-\rho_0^2)^2E_{20}}+J_D^\dagger/
\sqrt{m}+O_P(1/m),
\]
and an equation equivalent to (\ref{W+-EW+}) holds for $W^{\dag
+}-E^\dag W^{\dag+}$.

For some $0<c<C<\infty$, let
%
\begin{equation}
\label{calE} {\cal{E}}=\Biggl\{\mathbf{C}\dvtx \frac{1}{m+1}\sum
_{i=0}^mX_{2i+1}^8<C,
\frac
{1}{m+1}\sum_{i=0}^mX_{2i+1}^2>c
\Biggr\}.
\end{equation}
In Theorem~\ref{rne0condthm}, the $O_P(1/m)$, can be replaced by $\theta M_m$, where
$|\theta|<C$ and
\[
M_m=m\sum_{i=1}^mEY_i^4\Big/
\Biggl[\sum_{i=1}^mEY_i^2
\Biggr]^2
\]
as shown in Section 1 of the supplementary material of \citet{FRS13},
and for $\mathbf{C}\in{\cal{E}}$, $M_m$ is bounded. So
\[
P(R>u|{\cal{E}})=E\bigl[P(S>0|\mathbf{C})|{\cal{E}}\bigr]=E\bigl[\bar\Phi\bigl(
\sqrt {m}W^+\bigr)|{\cal{E}}\bigr]\bigl(1+O_P(1/m)\bigr).
\]
Using this and the equivalent term for $P^\dag(R^\dag>u|{\cal{E}})$, we have
%
\begin{eqnarray}
\label{ratR}&& \frac{|P(R>u|{\cal{E}})-P^\dag(R^\dag>u|{\cal{E}})|}{P(R>u|{\cal
{E}})}\nonumber\\
&&\qquad=\frac{|E^\dag[\bar\Phi(\sqrt{m}W^{\dag+})|{\cal{E}}]- E[\bar
\Phi(\sqrt{m}W^+)|{\cal{E}}]|}{ E[\bar\Phi(\sqrt{m}W^+)|{\cal
{E}}]}
\\
&&\qquad\leq\frac{I_1+I_2+I_3}{\bar\Phi(\sqrt{m}E(W^+|{\cal{E}}))},\nonumber
\end{eqnarray}
where we have used Jensen's inequality in the denominator and
%
\begin{eqnarray}
\label{I1} I_1&=&\bigl|\bar\Phi\bigl(\sqrt{m}E^\dag
\bigl(W^{\dag+}|{\cal{E}}\bigr)\bigr)- \bar\Phi\bigl(\sqrt {m}E\bigl(W^+|{
\cal{E}}\bigr)\bigr)\bigr|,
\\
\label{I2} I_2&=&\bigl|E^\dag\bigl[\bar\Phi\bigl(
\sqrt{m}W^{\dag+}\bigr)|{\cal{E}}\bigr]- \bar\Phi\bigl(\sqrt
{m}E^\dag\bigl(W^{\dag+}|{\cal{E}}\bigr)\bigr)\bigr|
\end{eqnarray}
and
%
\begin{equation}
\label{I3} I_3=\bigl|E\bigl[\bar\Phi\bigl(\sqrt{m}W^{+}
\bigr)|{\cal{E}}\bigr]- \bar\Phi\bigl(\sqrt {m}E\bigl(W^+|{\cal{E}}\bigr)\bigr)\bigr|.
\end{equation}

Noting that, for $\varphi(x)=-\bar\Phi'(x)$, $\varphi'(x)=-x\varphi
(x)$ and $x<\varphi(x)/\bar\Phi(x)<1+x$, we have
\[
\bar\Phi\bigl(\sqrt{m}E\bigl(W^+|{\cal{E}}\bigr)\bigr)>\varphi\bigl(\sqrt{m}E
\bigl(W^+|{\cal {E}}\bigr)\bigr)/\bigl(1+\sqrt{m}E\bigl(W^+|{\cal{E}}\bigr)\bigr).
\]
Then
\begin{eqnarray*}
\frac{I_3}{\bar\Phi(\sqrt{m}E(W^+|{\cal{E}}))}&\leq&\frac{m}{2}\frac
{E[(W^+-E(W^+|{\cal{E}}))^2\varphi(\sqrt{m}W^\ddag)|{\cal{E}}]}{\varphi
(\sqrt{m}E(W^+|{\cal{E}}))/(1+\sqrt{m}E(W^+|{\cal{E}}))},
\end{eqnarray*}
where $W^\ddag$ lies between $W^+$ and $E(W^+|{\cal{E}})$. Now, for
$\mathbf{C}\in{\cal{E}}$, noting (\ref{Wexp}) and (\ref{u0exp}),
\[
\frac{\varphi(\sqrt{m}W^\ddag)}{\varphi(\sqrt{m}E(W^+|{\cal
{E}}))}=O_P\bigl(e^{\sqrt{m}(u-\rho_0)^2}\bigr)=O_P(1)
\]
for $u=O(m^{-1/4})$, and using (\ref{Wexp}) and (\ref{W+-EW+}), we have
\[
\frac{I_3}{\bar\Phi(\sqrt{m}E(W^+|{\cal{E}}))}=O_P\bigl(m(u-\rho_0)^4+1/m
\bigr).
\]
An equivalent result holds for $I_2$. Also, using the same results gives
\begin{eqnarray*}
\frac{I_1}{\bar\Phi(\sqrt{m}E(W^+|{\cal{E}}))}&=&\frac{\sqrt{m}|E^\dag
(W^{\dag+}|{\cal{E}})-E(W^+|{\cal{E}})|\varphi(\sqrt{m}W^*)}{\varphi
(\sqrt{m}E(W^+|{\cal{E}}))/(1+\sqrt{m}E(W^+|{\cal{E}}))}
\\
&=&O_P\bigl(\sqrt{m}(u-\rho_0)^3+1/m
\bigr),
\end{eqnarray*}
where $W^*$ lies between $E(W^+|{\cal{E}})$ and $E^\dag(W^{\dag+}|{\cal{E}})$.

Finally, we need to consider the relative errors of the bootstrap and
the smoothed bootstrap.
%
\begin{Lemma}\label{le2}
For $\tau=O(1/\sqrt{m})$ and $u-\rho_0=O(n^{-1/4}),$
\[
P^\dagger\bigl(R^\dagger\geq u|{\cal{E}}\bigr)/P^*\bigl(R^*\geq
u|{\cal {E}}\bigr)=1+O_P\bigl(m(u-\rho_0)^4+1/m
\bigr).\vadjust{\goodbreak}
\]
\end{Lemma}

The proof of Lemma~\ref{le2} is given in Section 2 of the supplementary
material of \citet{FRS13}. Thus we have the following theorem:
%
\begin{thmm}\label{th2}
For ${\cal{E}}$ defined in (\ref{calE}), $u\geq\rho_0$, $u-\rho
_0=O(m^{-1/4})$ and $1-2u^2+\rho_0^2>\delta>0$,
\[
\frac{P(R>u|{\cal{E}})-P^*(R^*>u|{\cal{E}})}{P(R>u|{\cal
{E}})}=O_P\bigl(m(u-\rho_0)^4+1/m
\bigr).
\]
Further, if $E\varepsilon_1^8$ exists, then $P({\cal{E}})=1-o(1)$, if
$E\varepsilon_1^{16}$ exists, then $P({\cal{E}})=1-O(1/m)$ and if $\varepsilon
_1$ is bounded, then $P({\cal{E}})=1$, in which case the conditional
probabilities can be replaced by their expectations over ${\cal{E}}$.
\end{thmm}

\subsection{The conditional bootstrap}\label{sec4.3}

Consider obtaining a smoothed conditional bootstrap given $\mathbf{C}$. Let
\[
f_n(z)=\frac{1}{n-1}\sum_{k=2}^n
\frac{e^{-(z-\varepsilon_k)^2/2\tau
^2}}{\sqrt{2\pi\tau^2}},
\]
where $\varepsilon_i=X_i-\rho_0X_{i-1}$, for $i=2,\ldots,n$. Note that
this differs from $f_n$ of (\ref{fn}) in that the unstandardized errors
are used. Then the conditional density of $X_{2i}^\#$, the smoothed
bootstrap values of the even subscripted variable, given $X_{2i-1}$ and
$X_{2i+1}$ is

\[
g^\#(z|X_{2i-1},X_{2i+1})=\frac{f_n(z-\rho_0 X_{2i-1})f_n(X_{2i+1}-\rho
_0 z)}{\int f_n(z-\rho_0 X_{2i-1})f_n(X_{2i+1}-\rho_0 z)\,dz},
\]
where
\[
g^\#(z|X_{2i-1},X_{2i+1})=\frac{1}{(n-1)^2}\sum
_k\sum_l g_{ikl}^
\#(z),
\]
and, as in Section~\ref{appcond}, this can be reduced to
\[
g_{ikl}^\#(z)=\bigl(2\pi\tau^2/\bigl(1+
\rho_0^2\bigr)\bigr)^{-1/2}e^{-(1+\rho
_0^2)(z'-(\varepsilon_k-\rho_0\varepsilon_l)/(1+\rho_0^2))^2/2\tau
^2}w_{ikl}^
\#,
\]
where
\[
w_{ikl}^\#=\frac{e^{-(X_{2i+1}-\rho_0^2 X_{2i-1}-\rho_0\varepsilon
_k-\varepsilon_l)^2/2\tau^2(1+\rho_0^2)}} {
\sum_k\sum_l e^{-(X_{2i+1}-\rho_0^2 X_{2i-1}-\rho_0\varepsilon_k-\varepsilon
_l)^2/2\tau^2(1+\rho_0^2)}}
\]
and $z'=z-\rho_0(X_{2i-1}+X_{2i+1})/(1+\rho_0^2)$.

For each $i$ we sample from this distribution by first choosing
$\varepsilon_k,\varepsilon_l$ with probabilities $w_{ikl}^\#$, then obtaining
a random normal variable $Z_i'$ with mean $(\varepsilon_k-\rho_0\varepsilon
_l)/(1+\rho_0^2)$ and variance $\tau^2/(1+\rho_0^2)$, then taking
$X_{2i}^\#=Z_i'+\rho_0(X_{2i-1}+X_{2i+1})/(1+\rho_0^2)$.

Then the conditional cumulant generating function of $S^\#$ given
$\mathbf{C}$ is
\begin{eqnarray*}
mK^\#(t,u)&=&\sum_{i=1}^{m}\log \int
e^{\{t(A_i z-u(z^2+
B_i))\}}g^\#(z|X_{2i-1},X_{2i+1})\,dz
\\
&=&\sum_{i=1}^{m}\log \int
e^{-tu(z-A_i/2u)^2}g^\# (z|X_{2i-1},X_{2i+1})\,dz
 + m\frac{t(\bar{A}^{ 2}-4u^2\bar{B})}{4u}.
\end{eqnarray*}
Proceeding as in Section~\ref{appcond} we have
\[
K_{10}^\#(0,u)=\frac{1}{m} \sum_{i=1}^{m}
\int\bigl(zA_i-uz^2\bigr)g^\# (z|X_{2i-1},X_{2i+1})\,dz-u
\bar{B}.
\]
Let $u_{0}^\#$ be such that $K_{10}^\#(0,u_{0}^\#)=0$, then
%
\begin{equation}
u_0^\#=\frac{ \sum_{i=1}^{m}\int zg^\#(z|X_{2i-1},X_{2i+1})\,dz A_i}{\sum_{i=1}^{m}\int z^2g^\#(z|X_{2i-1},X_{2i+1})\,dz+m\bar{B}}.\label{u0th}
\end{equation}
So for $u> u_0^\#$,
\[
K_{10}^\#(0,u)=\bigl(u_0^\#-u\bigr) \Biggl(
\frac{1}{m}\sum_{i=1}^{m}\int
z^2g(z|X_{2i-1},X_{2i+1})\,dz+\bar{B}\Biggr)<0
\]
and $K_{20}^\#(t,u)>0$. So for $u>u_0^\#$, $K_{10}^\#(t,u)$ is
increasing in $t$, is negative for $t=0$ and as $t\to\infty$,
\[
K_{10}^\#(t,u)\to\frac{\bar{A}^{2}-4u^2\bar{B}}{4u}.
\]
Thus the saddle-point equation $K_{10}^\#(t,u)=0$ has a finite solution
$t^\#(u)$ for $u>u_0^\#$, if and only if $\bar{A}^{2}-4u^2\bar{B}>0$.
Further, $K^\#(t^\#(u),u)$ exists and is finite if $\bar{A}^{2}-4u^2\bar
{B}>0$. If $\bar{A}^{2}-4u^2\bar{B}<0$, $K^\#(t,u)\to-\infty$ as $t\to
\infty$.

Let $W^\#$, $W^{\#+}$ be defined in the same way as in the statement
of Theorem~\ref{rne0condthm}, then
%
\begin{equation}
P^\#\bigl(R^\#>u\bigr)=\bar\Phi\bigl(\sqrt{m}W^{\# +}\bigr)
\bigl(1+O_P(1/m)\bigr).\label{Sth1}
\end{equation}
Now
$K_{10}^\#(0,u_0^\#)=0$, so $t^{\#}(u_0^\#)=0$ and
\[
t^\#(u)=t^{\#\prime}\bigl(u_0^\#\bigr) \bigl(u-u_0^
\#\bigr)+\tfrac{1}{2}t^{\#\prime\prime
}\bigl(u_0^\#\bigr)
\bigl(u-u_0^\#\bigr)^2+O_P\bigl(
\bigl(u-u_0^\#\bigr)^3\bigr),
\]
with
\[
t^{\#\prime}\bigl(u_0^\#\bigr)=-K_{11}^
\#/K_{20}^\#,
\]
where $K_{ij}^\#=K_{ij}^\#(0,u_0).$ Then
\[
K_{11}^\#=\frac{1}{m} \sum_{i=1}^{m}
\int z^2g^\# (z|X_{2i-1},X_{2i+1})\,dz-\bar{B}
\]
and
\begin{eqnarray*}
K_{20}^\#&=&\frac{1}{m}\sum_{i=1}^m
\biggl\{\int\bigl(u_0^\# z^2-A_i z
\bigr)^2g^\# (z|X_{2i-1},X_{2i+1})\,dz
\nonumber
\\[-8pt]
\\[-8pt]
\nonumber
&&\hspace*{32pt}{} -\biggl[\int\bigl(u_0^\# z^2-A_i z
\bigr)g^\#(z|X_{2i-1},X_{2i+1})\,dz\biggr]^2\biggr\}.
\end{eqnarray*}

Now, as before, $D_1^\#=K_{11}^{\# 2}/2K_{20}^\#$. To compare $D_1^\#$
and $D_1$ we need the following lemma, the proof of which is given in Section 2
of the supplementary material \citet{FRS13}.

\begin{Lemma}\label{le3}
\[
\int h(z)g^\#(z|X_1,X_3)\,dz=\int h(z)g(z|X_1,X_3)\,dz
+O_P\biggl(\frac{1}{m}\biggr).
\]
\end{Lemma}

So, applying the lemma to $u_0^\#$, $K_{11}^\#$ and $K_{20}^\#$,
\[
D_1^\#=D_1+O_P(1/m).
\]

Now using (\ref{W+}) and an analogous term for $W^\#$ and noting that
$D_2-D^\#_2=O_P(1/\sqrt{m})$, we have
\[
\sqrt{m}\bigl(W^+-W^{\#+}\bigr) \bigl(1+\sqrt{m}W^+\bigr)=O\bigl(
\sqrt{m}(u-\rho_0)^3+1/m\bigr).
\]
Summarizing these results we have the following theorem:
%
\begin{thmm}\label{th3}
For $u\geq0$, $P(S>0|\mathbf{C})=0$ and $P(S^\#>0|\mathbf{C})=0$ if
$\bar{A}^{2}-4u^2\bar B<0$ and if $\bar{A}^{2}-4u^2\bar B>0$
$t(u)$ and $t^\#(u)$, solutions of $K_{10}(t,u)=0$ and $K^\#_{10}(t,u)=0$,
exist and are both finite and positive, and if $EX_1^8$ is bounded,
(\ref{Sth1}) holds and
\[
P(R>u|\mathbf{C})=P^\#\bigl(R^\#>u|\mathbf{C}\bigr)\bigl[1+O_P
\bigl(\sqrt{m}(u-\rho_0)^3+1/m\bigr)\bigr].
\]
\end{thmm}

\section{Numerical results}\label{sec5}

Monte Carlo simulations, bootstraps and tail area
approximations both unconditionally and conditionally are used to
illustrate accuracy of results
and to compare the power of the unconditional and the conditional bootstrap.

First we describe the computational methods.
The true distribution of $\hat\rho$ is approximated by Monte
Carlo simulations of 1,000,000.
For the bootstrap, we consider
testing $H_0\dvtx \rho=\rho_0$. The unconditional bootstrap is straightforward
in that we compute $n-1$ residuals, $\varepsilon_i=x_i-\rho_0 x_{i-1}$,
center them and
sample these with replacement. Then $x_i^*=\rho_0 x_{i-1}^*
+\varepsilon_i^*$ with $x_1^*=\varepsilon_1^*/(1-\rho^2)$, and we compute $R^*$
and obtain an estimate of $P^*(R^*>u)$ from repetitions. For the
conditional bootstrap
of Section~\ref{relerrbsnonzero},
we draw samples
$\varepsilon^{\dagger}_i$'s from $f_n$ in (\ref{fn}) with $\tau$ equal to
$1/m$. We first generate $X_i^{\dagger}$'s from the
$\varepsilon_i^{\dagger}$'s. Then $X_{2i}^{\dagger}$ are replaced by
generating an observation from the normal mixture given in
(\ref{gdag})--(\ref{gikldag}), $R^{\dagger}$ is computed and
repetitions give an estimate of $P^\dagger(R^\dagger>u|\mathbf{C^*})$.
Now repeating this entire process from sampling $\varepsilon^\dagger_i$'s
and averaging the conditional probabilities gives an estimate of
$P^\dagger(R^\dagger>u)$. For the conditional bootstrap of Section~\ref{sec4.3},
we replace $X_{2i}$ by $X_{2i}^\#$ drawn from~(\ref{gdag}), calculate
$R^\#$ and repeat this process to get an estimate of $P^\#(R^\#
>u|\mathbf{C})$.

\begin{table}
\caption{Comparison of saddle-point and simulated tail areas for normal
distribution from Section~\protect\ref{sec3} with the unconditional case (U) and the
conditional case (C) at both $n=39$ and $n=9$}\label{table1}
\begin{tabular*}{\textwidth}{@{\extracolsep{\fill}}lcccccc@{}}
\hline
&&\multicolumn{5}{c@{}}{\textbf{Tail prob. exceeds}}\\[-6pt]
&&\multicolumn{5}{c@{}}{\hrulefill}\\
\multicolumn{1}{@{}l}{}&  & $\bolds{\rho+0.05}$& $ \bolds{\rho+0.10}$&$\bolds{ \rho+0.15}$&
$\bolds{\rho+0.20}$& $\bolds{\rho+0.25}$\\
\hline
$n=39$&$\rho=0.5$&&&&\\
UC&saddle-point&0.3210&0.1923& 0.0946&0.0353&0.0088\\
&simulations&0.3223&0.1922& 0.0946& 0.0352&0.0094\\[3pt]
$n=9$&$\rho=0.5$&&&&\\
UC&saddle-point&0.3629&0.2937&0.2261&0.1624&0.1066\\
&simulations&0.3695&0.2994&0.2310&0.1660&0.1081\\[6pt]
$n=39$&$\rho=0.5$&&&&\\
C&saddle-point&0.3133&0.1888&0.0983&0.0412&0.0118\\
&simulation&0.3136&0.1884&0.0983&0.0410&0.0118\\[3pt]
$n=9$&$\rho=0.5$&&&&\\
C&saddle-point&0.4077& 0.3413& 0.2713& 0.1972& 0.1267\\
&simulation&0.4094& 0.3432& 0.2722& 0.1999& 0.1280\\
\hline
\end{tabular*}      \vspace*{-3pt}
\end{table}

The results for the approximations of Section~\ref{sec3} for the Gaussian case are
given in the upper part of Table~\ref{table1} for the unconditional
results (U) and the lower part for the conditional case (C). As can be seen,
the agreements between the simulation results and the
saddle-point, computed as in Section~\ref{sec3} for normal data, are
excellent with very accurate results, even for $n=9$. The accuracy for
values of
$\rho<0.5$ is even better.

\begin{table}
\caption{Unconditional bootstrap (BS: 100,000 replicates) and expected
conditional bootstrap averages over $\mathbf{C}^\dagger$ (ECBS: using
500 sets of the conditional bootstrap with 10,000 replicates) and
average of conditional saddle-point approximation (ECSP: over 500
replicates), from the same original sample from $t_{10}$}
\label{table2}
\begin{tabular*}{\textwidth}{@{\extracolsep{\fill}}lcccccc@{}}
\hline
&&\multicolumn{5}{c@{}}{\textbf{Tail prob. exceeds}}\\[-6pt]
&&\multicolumn{5}{c@{}}{\hrulefill}\\
\multicolumn{1}{@{}l}{}&  & $\bolds{\rho+0.05}$& $ \bolds{\rho+0.10}$&$\bolds{ \rho+0.15}$&
$\bolds{\rho+0.20}$& $\bolds{\rho+0.25}$\\
\hline
$n=39$&$\rho=0.5$&&&&\\
&BS& 0.3206&0.1921& 0.0943 &0.0350 &0.0086\\
&ECBS&0.3160&0.1833&0.0860&0.0309&0.0075 \\
&ECSP&0.3131&0.1823&0.0861&0.0308&0.0075\\
\hline
\end{tabular*}
\end{table}

In Table~\ref{table2}, we use a single sample from a $t_{10}$ distribution to
compare the unconditional bootstrap and the smoothed
bootstrap averaged over $\mathbf{C}^{\dagger}$'s for $\rho_0=0.5$,
to demonstrate the results of Lemma~\ref{le2}, and we
obtain an estimate of $E^\dagger_+\bar\Phi(\sqrt{m}W^{\dagger+})$, the
expected value of the saddle-point approximation given in Theorem~\ref{rne0condthm}, by
averaging over 100 values of $\mathbf{C}^\dagger$, comparing this to
the Monte Carlo estimates. These results, which would vary from sample
to sample from the $t_{10}$ distribution, illustrate excellent relative
accuracy, and we note that better results are obtained for $0\leq\rho_0<0.5$.

\begin{table}[b]
\caption{Simulated tail probabilities (SIM: 1,000,000 samples), estimates
of expected bootstrap tail probabilities and standard deviations of
bootstrap tail probabilities based on means and standard deviations of
40 samples (EBS and SDBS: 100,000 bootstrap replications) from $t_{10}$
and centered exponential distributions}
\label{table3}
\begin{tabular*}{\textwidth}{@{\extracolsep{\fill}}lcccccc@{}}
\hline
&&\multicolumn{5}{c@{}}{\textbf{Tail prob. exceeds}}\\[-6pt]
&&\multicolumn{5}{c@{}}{\hrulefill}\\
\multicolumn{1}{@{}l}{}&  & $\bolds{\rho+0.05}$& $ \bolds{\rho+0.10}$&$\bolds{ \rho+0.15}$&
$\bolds{\rho+0.20}$& $\bolds{\rho+0.25}$\\
\hline
$n=39$&$\rho=0.5$&&&&\\
$t_{10}$&SIM&0.3171&0.1885&0.0916&0.0340&0.0083\\
&EBS&0.3215&0.1932&0.0957&0.0361&0.0094 \\
&SDBS&0.0016&0.0017&0.0019&0.0015&0.0009\\[3pt]
exp&SIM&0.3174&0.1937&0.1020&0.0442&0.0154\\
&EBS&0.3223& 0.1991& 0.1059& 0.0473& 0.0173 \\
&SDBS&0.0044&0.0088&0.0123&0.0121&0.0089\\
\hline
\end{tabular*}\vspace*{-3pt}
\end{table}

In Table~\ref{table3}, to illustrate the main results of Theorem~\ref{th2}, we
compare the simulated distribution, when sampling from the
$t_{10}$-distribution and\vadjust{\goodbreak}
the exponential distribution shifted to have mean 0, with the bootstrap
averages over 40 samples. The average bootstrap is quite accurate,
while the standard deviation shows that the relative error of the
bootstrap becomes larger in the tails, as expected since this is shown
to be of order $m(u-\rho_0)^4$ in Theorem~\ref{th2}. For $0\leq\rho
_0<0.5$, there is even better accuracy.

Table~\ref{table4} illustrates the accuracy of the results of Theorem~\ref{th3}
using random samples for $\rho_0$ equal to 0 and 0.5 for centered
exponential errors. The saddle-point approximation has the relative
accuracy property. In this case, there is considerable variation in
tail areas as different random samples are taken, but similar accuracy
is achieved with other samples. Similar results are obtained for the
$t_{10}$ distribution and for $0\leq\rho_0<0.5$.

\begin{table}
\caption{
Comparison of tail areas for conditional bootstrap (CBS) and\break conditional saddle-point tail area (CSP)
for one sample from a centered exponential with $\rho_0=0.0$ and   $\rho_0=0.5$, as in Section \protect\ref{sec4.3}}
\label{table4}
\begin{tabular*}{\textwidth}{@{\extracolsep{\fill}}lcccccc@{}}
\hline
&&\multicolumn{5}{c@{}}{\textbf{Tail prob. exceeds}}\\[-6pt]
&&\multicolumn{5}{c@{}}{\hrulefill}\\
\multicolumn{1}{@{}l}{}& $\bolds{\rho}$ & $\bolds{\rho+0.05}$& $ \bolds{\rho+0.10}$&$\bolds{ \rho+0.15}$&
$\bolds{\rho+0.20}$& $\bolds{\rho+0.25}$\\
\hline
\multicolumn{1}{@{}l}{${n=39}$}&&&&&&\\
CSP&0.0&0.4300& 0.3103& 0.2074& 0.1268& 0.0697\\
CBS&0.0&0.4378& 0.3147& 0.2080& 0.1274& 0.0688\\[3pt]
CSP&0.5&0.2499& 0.0863& 0.0145& 0.0004 & 0.0000\\
CBS&0.5&0.2456& 0.0843& 0.0132& 0.0002 & 0.0000\\
\hline
\end{tabular*}\vspace*{-3pt}
\end{table}

\begin{table}
\caption{Power under unconditional (U) and conditional (C) tests for
the Gaussian case in the left half of the table and the general case
from $t_{10}$ in the right half}
\label{table5}
\begin{tabular*}{\textwidth}{@{\extracolsep{\fill}}l@{\hspace*{40pt}}cccc@{\hspace*{40pt}}cccc@{}}
\hline
&\textbf{U}&\textbf{C}&\textbf{U}&\textbf{C}&\textbf{U}&\textbf{C}&\textbf{U}&\textbf{C}\\
\multicolumn{1}{r}{$\bolds{\rho_0}=$}& \textbf{0}& \textbf{0} & \textbf{0.4}&\textbf{0.4}&\textbf{0}& \textbf{0} & \textbf{0.4}&\textbf{0.4}\\
\hline
$\rho_1=\rho_0+0.1$& 0.15& 0.15& 0.18 &0.12&0.15&0.13&0.18&0.11\\
$\rho_1=\rho_0+0.3$& 0.58 &0.59& 0.73& 0.42&0.58&0.53&0.73&0.38\\
$\rho_1=\rho_0+0.5$& 0.92 &0.90& 0.98& 0.89&0.93&0.90&0.98&0.78\\
\hline
\end{tabular*}
\end{table}

Finally, we compare the power of the two tests based on the
unconditional bootstrap and the conditional bootstrap in Table~\ref{table5} for
the Gaussian case of Section~\ref{sec3} and for the general case from Sections
\ref{relerrbsnonzero} and~\ref{sec4.3}. We note that the tests have equal power up to
computational accuracy when $\rho_0=0$, as might be expected since
there is no loss of information due to conditioning in this case, but
there is some loss of power in the case of $\rho_0=0.2$ and a
considerable loss for $\rho_0=0.4$.


\begin{supplement}[id=suppA]
\stitle{Supplement to ``Relative errors for bootstrap approximations of the
serial correlation coefficient''}
\slink[doi]{10.1214/13-AOS1111SUPP} 
\sdatatype{.pdf}
\sfilename{aos1111\_supp.pdf}
\sdescription{We provide details and proofs needed for a number of results in the paper.}
\end{supplement}

%

\printaddresses


\begin{thebibliography}{8}

\bibitem[\protect\citeauthoryear{Anderson}{1959}]{A59}
\begin{barticle}[mr]
\bauthor{\bsnm{Anderson},~\bfnm{T.~W.}\binits{T.~W.}}
(\byear{1959}).
\btitle{On asymptotic distributions of estimates of parameters of stochastic
  difference equations}.
\bjournal{Ann. Math. Statist.}
\bvolume{30}
\bpages{676--687}.
\bid{issn={0003-4851}, mr={0107347}}
\bptok{imsref}%
\end{barticle}
\endbibitem

\bibitem[\protect\citeauthoryear{Bose}{1988}]{BA}
\begin{barticle}[mr]
\bauthor{\bsnm{Bose},~\bfnm{Arup}\binits{A.}}
(\byear{1988}).
\btitle{Edgeworth correction by bootstrap in autoregressions}.
\bjournal{Ann. Statist.}
\bvolume{16}
\bpages{1709--1722}.
\bid{doi={10.1214/aos/1176351063}, issn={0090-5364}, mr={0964948}}
\bptok{imsref}%
\end{barticle}
\endbibitem

\bibitem[\protect\citeauthoryear{Daniels}{1956}]{D56}
\begin{barticle}[mr]
\bauthor{\bsnm{Daniels},~\bfnm{H.~E.}\binits{H.~E.}}
(\byear{1956}).
\btitle{The approximate distribution of serial correlation coefficients}.
\bjournal{Biometrika}
\bvolume{43}
\bpages{169--185}.
\bid{issn={0006-3444}, mr={0079395}}
\bptok{imsref}%
\end{barticle}
\endbibitem

\bibitem[\protect\citeauthoryear{Field and Robinson}{2013}]{FRS13}
\begin{bmisc}[author]
\bauthor{\bsnm{Field},~\bfnm{Chris}\binits{C.}} \AND
  \bauthor{\bsnm{Robinson},~\bfnm{John}\binits{J.}}
(\byear{2013}).
\bhowpublished{Supplement to ``Relative errors for bootstrap approximations of
  the serial correlation coefficient.'' DOI:\doiurl{10.1214/13-AOS1111SUPP}.}
\bptok{imsref}%
\end{bmisc}
\endbibitem

\bibitem[\protect\citeauthoryear{Lieberman}{1994a}]{L89}
\begin{barticle}[mr]
\bauthor{\bsnm{Lieberman},~\bfnm{Offer}\binits{O.}}
(\byear{1994}a).
\btitle{Saddlepoint approximation for the distribution of a ratio of quadratic
  forms in normal variables}.
\bjournal{J. Amer. Statist. Assoc.}
\bvolume{89}
\bpages{924--928}.
\bid{issn={0162-1459}, mr={1294736}}
\bptok{imsref}%
\end{barticle}\vadjust{\goodbreak}
\endbibitem

\bibitem[\protect\citeauthoryear{Lieberman}{1994b}]{L94}
\begin{barticle}[mr]
\bauthor{\bsnm{Lieberman},~\bfnm{Offer}\binits{O.}}
(\byear{1994}b).
\btitle{Saddlepoint approximation for the least squares estimator in
  first-order autoregression}.
\bjournal{Biometrika}
\bvolume{81}
\bpages{807--811}.
\bid{doi={10.1093/biomet/81.4.807}, issn={0006-3444}, mr={1326430}}
\bptok{imsref}%
\end{barticle}
\endbibitem

\bibitem[\protect\citeauthoryear{Phillips}{1978}]{P78}
\begin{barticle}[author]
\bauthor{\bsnm{Phillips},~\bfnm{P.~C.~B.}\binits{P.~C.~B.}}
(\byear{1978}).
\btitle{Edgeworth and saddle-point approximations in the first-order
  noncircular autoregression}.
\bjournal{Biometrika}
\bvolume{65}
\bpages{91--98}.
\bptok{imsref}%
\end{barticle}
\endbibitem

\end{thebibliography}
\end{document}